\documentclass{article}
\usepackage{graphics, graphicx, latexsym, amssymb}
\newtheorem{theo}{Theorem}
\newtheorem{coro}{Corollary}
\newtheorem{proc}{Procedure}

\begin{document}

\title{The area generating function \\ for simplex-duplex polyominoes}
\author{Svjetlan Fereti\'{c} \footnote{e-mail: svjetlan.feretic@gradri.hr}\\ Faculty of Civil Engineering, University of Rijeka, \\ Viktora Cara Emina 5, 51000 Rijeka, Croatia}
\maketitle

\begin{abstract}

Back in the early days of polyomino enumeration, a model called column-convex polyominoes was introduced and its area generating function was found. That generating function is rational: the numerator has degree four and the denominator has degree three. Let a column-duplex polyomino be a polyomino whose columns can have either one or two connected components. In spite of being an immediate generalization of column-convex polyominoes, column-duplex polyominoes do not seem amenable to any kind of exact enumeration. Hence, we are interested in models which interpolate between column-convex polyominoes and column-duplex polyominoes. This paper is concerned with a model which we have called simplex-duplex polyominoes. A simplex-duplex polyomino is a column-duplex polyomino in which there is no occurrence of two adjacent columns each having two connected components. Simplex-duplex polyominoes are not easy to deal with, but their area generating fuction can still be found. To find this generating function, we use an upgraded version of the Temperley method. Though that technique is widely used in these times, our application presents two interesting features. Firstly, we add one or two columns at a time, thus bypassing those simplex-duplex polyominoes which end with a two-component column. (It is somewhat more usual to add just one column at a time.) Secondly, we obtain a functional equation that involves both the first and the second derivatives of the sought-for generating function. (Such equations usually involve the first derivative only. In some cases, no derivative is involved at all.) Right because of this latter interesting feature, the Temperley method produces a very complicated formula for the generating function. Anyway, from that formula it is easy to compute Taylor polynomials. Thus we get plenty of evidence that the number of $n$-celled simplex-duplex polyominoes behaves asymptotically as $0.119443 \cdot 3.522020^n$. For comparison, the number of $n$-celled column-convex polyominoes behaves asymptotically as $0.180916 \cdot 3.205569^n$.

\end{abstract}

\section{Introduction}

Enumeration of polyominoes is a well-established topic in mathematics, physics and chemistry \cite{book}. It would be ideal to know how many polyominoes can be made up of $n$ cells, or how many borders of polyominoes can be made up of $m$ unit edges. Realistically speaking, however, neither of these questions will be answered in the foreseeable future. So it is a sound idea to confine attention to sets of polyominoes which do not contain all polyominoes, but anyway contain rather many of them. Sets of that kind sometimes allow of the enumeration by perimeter and/or by area. Yet, in defining such a ``friendly" set of polyominoes, there is not much room left for poetic license. That kind of definitions almost always require that the elements of the set be either in a certain sense \textit{convex} or in a certain sense \textit{directed}. A classical example are \textit{column-convex polyominoes}: their definition requires that the intersection between the polyomino and any vertical straight line must be a connected set. The area generating function for column-convex polyominoes has been known at least since 1956 \cite{Polya, Temperley}. The perimeter generating function for column-convex polyominoes is also not a recent result \cite{Brak, Delest, F_1996}. However, there exists an ongoing line of research aimed at constructing models which are supersets of column-convex polyominoes and are still solvable. The models constructed so far are called $m$\textit{-convex polygons} \cite{m-convex}, \textit{prudent polygons} \cite{prudent}, \textit{cheesy polyominoes} \cite{semi}, \mbox{\textit{polyominoes with cheesy blocks} \cite{semiblo},} and \textit{column-subconvex} \mbox{\textit{polyominoes} \cite{undi}.} The former two models allow of the enumeration by perimeter, whereas the latter three models admit of the enumeration by area. (Thus, it is still an open problem to think out a model which significantly generalizes column-convex polyominoes and admits of both the enumeration by perimeter and the enumeration by area.) 

Let a \textit{column-duplex polyomino} be a polyomino whose columns can have either one or two connected components. Somewhat surprisingly, column-duplex polyominoes are already too hard a model: they seem to admit of no exact enumeration. This paper is concerned with a model which we have called \textit{simplex-duplex polyominoes}. A cheesy polyomino, a polyomino with cheesy blocks, a column-subconvex polyomino and a simplex-duplex polyomino are all column-duplex polyominoes. In the first three models, the definition requires (perhaps among other things) that the gap within a two-component column must not be greater than $m$ cells in size, where $m$ is a positive integer which we fix in advance. In the fourth model (simplex-duplex polyominoes), the gap within a column can be of any size, but it is not allowed that two two-component columns come one immediately after the other.

The main result of this paper is a formula for a bivariate generating function $G(q,w)$, in which the coefficient of $q^n w^k$ is the number of simplex-duplex polyominoes having $n$ cells and $k$ two-component columns. To compute the formula for $G(q,w)$, we make use of an upgraded version of the Temperley method 
\cite{Bousquet, Temperley}. Our computations have some interesting features, as is described in the abstract above. Although $G(q,w)$ and $G(q,1)$ prove to be very complicated $q$-series, from the formulas for $G(q,w)$ and $G(q,1)$ it is easy to find Taylor polynomials of any reasonable degree. Thus we get plenty of evidence that the number of $n$-celled simplex-duplex polyominoes behaves asymptotically as $0.119443 \cdot 3.522020^n$. This means that simplex-duplex polyominoes are one of the largest sets of polyominoes for which the area generating function is known (cf. \cite{animals}). We already enumerated some fairly large sets of polyominoes \cite{semi, semiblo, undi}, but those polyominoes had hexagonal cells, whereas simplex-duplex polyominoes have square cells. However, it is fair to say that the best currently known lower bound on Klarner's constant is $3.980137$ \cite{twisted}, which is a significantly greater number than our growth constant $3.522020$.

This paper continues as follows. In Section 2, we state the necessary definitions and conventions. In Section 3, we establish and solve a functional equation. The function satisfying that equation is a generating function not for all simplex-duplex polyominoes, but only for those simplex-duplex polyominoes which end with a one-component column. Once the functional equation is solved, it is however very easy to obtain the area and duplex columns generating function for all simplex-duplex polyominoes. The formula for this latter generating function is stated in Section 4. In the same section, there is also a corollary, in which we state the area generating function for column-convex polyominoes. In Section 5, there is a bit of asymptotic analysis. In Section 6 we conclude, outlining further work prompted by our results.

\section{Definitions and conventions}

In this paper, a \textit{cell} is a unit square whose vertices have integer coordinates. A plane figure $P$ is a \textit{polyomino} if $P$ is a union of finitely many cells and the interior of $P$ is connected. 

Let $P$ and $Q$ be two polyominoes. We consider $P$ and $Q$ to be equal if and only if there exists a translation $f$ such that $f(P)=Q$.

Given a polyomino $P$, it is useful to partition the cells of $P$ according to their horizontal projection. Each block of that partition is a \textit{column} of $P$. With this definition, a column of a polyomino is not necessarily a connected set. However, it may happen that every column of a polyomino $P$ is a connected set. In this case, $P$ is called a \textit{column-convex} polyomino. 

By a \textit{column-duplex} polyomino, we mean a polyomino in which columns with three or more connected components are not allowed. Thus, each column of a column-duplex polyomino has either one or two connected components. 

Let $P$ be a column-duplex polyomino and let $c$ be a column of $P$. If $c$ has one connected component, we say that $c$ is a \textit{simplex} column. If $c$ has two connected components, we say that $c$ is a \textit{duplex} column.

A \textit{simplex-duplex} polyomino is such a column-duplex polyomino in which consecutive duplex columns are not allowed. If $c$ is a column of a simplex-duplex polyomino, and $c$ is a (whether left or right) neighbour of a duplex column, then $c$ must be a simplex column. See Figure 1.

\begin{figure}
\begin{center}
\includegraphics[width=63mm]{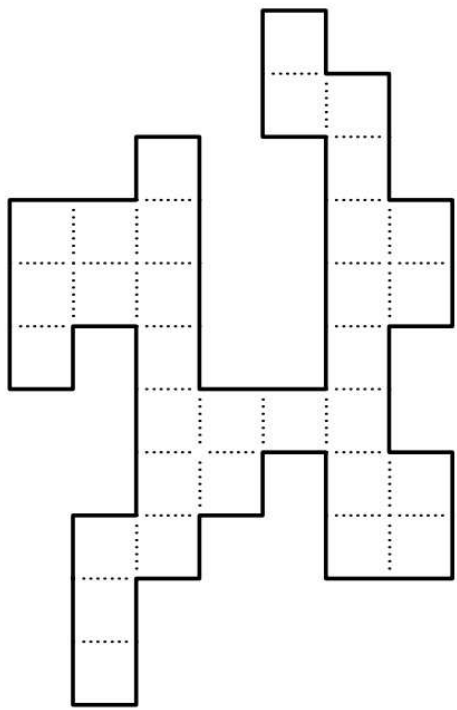} 
\caption{A simplex-duplex polyomino.}
\end{center}
\end{figure}

Let $P$ be a simplex-duplex polyomino and let $P$ have at least two simplex columns. Then we define the \textit{lower pivot cell} of $P$ to be the cell which is the right neighbour of the bottom cell of the second-last (\textit{i.e.}, second-rightmost) among the simplex columns of $P$. We also define the \textit{upper pivot cell} of $P$ to be the right neighbour of the top cell of the second-last among the simplex columns of $P$.

Let $P$ be a simplex-duplex polyomino and let $P$ have at least one duplex column. Then we define the \textit{lower inner pivot cell} of $P$ to be the right neighbour of the top cell of the lower component of the last among the duplex columns of $P$. We also define the \textit{upper inner pivot cell} of $P$ to be the right neighbour of the bottom  cell of the upper component of the last among the duplex columns of $P$. 

Observe that the lower pivot cell of a simplex-duplex polyomino $P$ is not necessarily contained in $P$. The same holds for the upper pivot cell, the lower inner pivot cell and the upper inner pivot cell of $P$. 

If a polyomino $P$ is made up of $n$ cells, we say that the \textit{area} of $P$ is $n$. 

Let $S$ denote the set of those simplex-duplex polyominoes whose last (\textit{i.e.}, rightmost) column is a simplex column. It is convenient to first compute a generating function for the set $S$, and thence a generating fuction for all simplex-duplex polyominoes. So, let

\begin{displaymath}
A(q,t,w)=\sum_{P \in S} q^{area \ of \ P} \cdot t^{{{the \ height \ of \ the \atop last \ column \ of \ P} \atop }} \cdot 
w^{{{the \ number \ of \atop duplex \ columns \ of \ P} \atop}}.
\end{displaymath}

Next, we define three generating functions in two variables, $q$ and $w$: Let $A_1=A(q,1,w)$, $B_1=\left\{\frac{\partial \left[\frac{A(q,t,w)}{t}\right]}{\partial t}\right\}_{with \ t=1}$ and $C_1=\frac{1}{2} \cdot \left\{\frac{\partial^2 \left[\frac{A(q,t,w)}{t}\right]}{\partial t^2}\right\}_{with \ t=1}$. \vspace{1mm}

For $n \in \mathbb{N}$, let $T_n$ be the set of simplex-duplex polyominoes having a duplex last column and an $n$-celled second-last column. Every element of $T_n$
is made up of three ingredients:

\begin{enumerate}
\item the ``left factor" consisting of all columns except the last column,
\item the lower connected component of the last column,
\item the upper connected component of the last column.
\end{enumerate}

The area and duplex columns generating function for ingredient No. 1 is $\langle t^n  \rangle A(q,t,w)$, where $\langle t^n  \rangle$ stands for ``the coefficient of $t^n$". The area generating function for each of the other two ingredients is $\frac{q}{1-q}$. The three ingredients can be glued together in exactly ${n-1 \choose 2}=\frac{(n-1)(n-2)}{2}$ ways, because there are ${n-1 \choose 2}$ ways to choose two non-adjacent cells in the last column of ingredient No. 1. Altogether, the area and duplex columns generating function for the set $T_n$ is $\left[ \langle t^n  \rangle A(q,t,w) \right] \cdot \frac{(n-1)(n-2)}{2} \cdot \left( \frac{q}{1-q} \right)^2 \cdot w$. As to all simplex-duplex polyominoes ending with a duplex column, their area and duplex columns generating function is

\begin{eqnarray*}
& \sum_{n=1}^{\infty} \left[ \langle t^n  \rangle A(q,t,w) \right] \cdot \frac{(n-1)(n-2)}{2} \cdot \left( \frac{q}{1-q} \right)^2 \cdot w \\
& = \frac{q^2 w}{(1-q)^2} \cdot \sum_{n=1}^{\infty} \frac{(n-1)(n-2)}{2} \cdot \left[ \langle t^n  \rangle A(q,t,w) \right] \\
& = \frac{q^2 w}{(1-q)^2} \cdot C_1.
\end{eqnarray*}

Let

\begin{displaymath}
G=G(q,w)=\sum_{{P \ a \ simplex-duplex \atop \ polyomino}} q^{area \ of \ P} \cdot w^{{{the \ number \ of \atop duplex \ columns \ of \ P} \atop}}.
\end{displaymath} 

By now it is clear that 

\begin{equation}
\label{ekg}
G=A_1 + \frac{q^2 w}{(1-q)^2} \cdot C_1.
\end{equation}

\section{Establishing and solving a functional equation}
 
Henceforth the notation $A(q,t,w)$ will be abbreviated as $A(t)$.

In order to obtain a functional equation for the generating function $A(t)$, we are going to suitably partition the set $S$. (Recall that $S$ is the set of those simplex-duplex polyominoes whose last column is a simplex column.)

First, let $S_\alpha$ be the set of those $P \in S$ which have no other simplex column than the last column.

Let $S_\beta$ be the set of those $P \in S \setminus S_{\alpha}$ which have the following two properties:

\begin{itemize}
\item the second-last column is a simplex column,
\item the last column contains the lower pivot cell of $P$.
\end{itemize}

Let $S_\gamma$ be the set of $P \in S\setminus S_{\alpha}$ having the following two properties:

\begin{itemize}
\item the second-last column is a simplex column,
\item the last column does not contain the lower pivot cell of $P$.
\end{itemize}

Thus, $S_{\beta} \cup S_{\gamma}$ is the set of those $P \in S \setminus S_{\alpha}$ whose second-last column is a simplex column.

Let $S_\delta$ be the set of those $P \in S \setminus S_{\alpha}$ which have the following three properties:

\begin{itemize}
\item the second-last column is a duplex column and the third-last column is (necessarily) a simplex column,
\item the lower component of the second-last column and the third-last column have no edge in common,
\item the lower pivot cell of $P$ is contained in the upper component of the second-last column.
\end{itemize}

Let $S_\epsilon$ be the set of $P \in S \setminus S_{\alpha}$ having the following three properties:

\begin{itemize}
\item the second-last column is a duplex column and the third-last column is a simplex column,
\item the lower component of the second-last column and the third-last column have no edge in common,
\item the lower pivot cell of $P$ is contained in the hole of the second-last column.
\end{itemize}

The definition of $S_\zeta$ is obtained from the definition of $S_\delta$ by writing the word ``upper" where the definition of $S_\delta$ says ``lower", and by writing the word ``lower" where the definition of $S_\delta$ says ``upper". 

The definition of $S_\eta$ is obtained when the changes just described are made to the definition of $S_\epsilon$ (instead of to the definition of $S_\delta$). However, since the set $S_\eta$ appears in the proof of Theorem 1 below, it makes sense to define $S_\eta$ explicitly: $S_\eta$ is the set of $P \in S \setminus S_{\alpha}$ having the following three properties:

\begin{itemize}
\item the second-last column is a duplex column and the third-last column is a simplex column,
\item the upper component of the second-last column and the third-last column have no edge in common,
\item the upper pivot cell of $P$ is contained in the hole of the second-last column.
\end{itemize}

Thus, $S_{\delta} \cup S_{\epsilon} \cup S_{\zeta} \cup S_{\eta}$ is the set of those $P \in S \setminus S_{\alpha}$ which, in addition to having a duplex second-last column, also have the property that $P \setminus (the \ last \ column \ of \ P)$ is not a polyomino. 

Let $S_\theta$ be the set of $P \in S \setminus S_{\alpha}$ which have the following three properties:

\begin{itemize}
\item the second-last column is a duplex column and the third-last column is a simplex column,
\item each of the two components of the second-last column has at least one edge in common with the third-last column,
\item the last column also has at least one edge in common with each of the two components of the second-last column.
\end{itemize}

Let $S_\iota$ be the set of $P \in S \setminus S_{\alpha}$ which have the following four properties:

\begin{itemize}
\item the second-last column is a duplex column and the third-last column is a simplex column,
\item each of the two components of the second-last column has at least one edge in common with the third-last column,
\item the last column has at least one edge in common with the lower component of the second-last column, but does not have any edges in common with the upper component of the second-last column,
\item the last column does not contain the lower inner pivot cell of $P$.
\end{itemize}

Let $S_\kappa$ be the set of $P \in S \setminus S_{\alpha}$ having the following four properties:

\begin{itemize}
\item the second-last column is a duplex column and the third-last column is a simplex column,
\item each of the two components of the second-last column has at least one edge in common with the third-last column,
\item the last column has at least one edge in common with the lower component of the second-last column, but does not have any edges in common with the upper component of the second-last column,
\item the last column contains the lower inner pivot cell of $P$.
\end{itemize}

The definition of $S_\lambda$ is obtained from the definition of $S_\iota$ by writing the word ``upper" where the definition of $S_\iota$ says ``lower", and by writing the word ``lower" where the definition of $S_\iota$ says ``upper". 

The definition of $S_\mu$ is obtained when the changes just described are made to the definition of $S_\kappa$ (instead of to the definition of $S_\iota$).

Thus, $S_{\theta} \cup S_{\iota} \cup S_{\kappa} \cup S_{\lambda} \cup S_{\mu}$ is the set of those $P \in S \setminus S_{\alpha}$ which, in addition to having a duplex second-last column, also have the property that $P \setminus (the \ last \ column \ of \ P)$ is a polyomino. This means that $S_{\delta} \cup S_{\epsilon} \cup 
\ldots \cup S_{\mu}$ is the set of all $P \in S \setminus S_{\alpha}$ whose second-last column is a duplex column. 

The sets $S_\alpha, \ S_\beta, \ldots, \ S_\mu$ form a partition of the set $S$.

Let the notations $A_\alpha(t), \ A_\beta(t), \ldots, \ A_\mu(t)$ denote the parts of the series $A(t)$ that come from the sets $S_\alpha, \ S_\beta, \ldots, \ S_\mu$, respectively.

\begin{theo} We have:
\begin{eqnarray*}
A_\alpha(t) & = & \frac{qt}{1-qt}+\frac{q^5t^3w}{(1-q)^2(1-qt)^3} \ , \\
A_\beta(t) & = & \frac{qt}{(1-qt)^2} \cdot A_1 \ , \\
A_\gamma(t) & = & \frac{qt}{1-qt} \cdot B_1 \ , \\
A_\delta(t)=A_\zeta(t) & = & \frac{q^5t^3w}{(1-q)^3(1-qt)^3} \cdot A_1 \ , \\
A_\epsilon(t)=A_\eta(t) & = & \frac{q^5t^3w}{(1-q)^2(1-qt)^4} \cdot A_1 - \frac{q^4t^2w}{(1-q)^2(1-qt)^4} \cdot A(qt) \ , \\
A_\theta(t) & = & \frac{q^5t^3w}{(1-q)^2(1-qt)^3} \cdot B_1 \\ 
& & \mbox{} - \frac{q^5t^3w}{(1-q)^2(1-qt)^4} \cdot A_1 + \frac{q^4t^2w}{(1-q)^2(1-qt)^4} \cdot A(qt) \ , \\
A_\iota(t)=A_\lambda(t) & = & \frac{q^4tw}{(1-q)^3(1-qt)} \cdot C_1 \ , \\
A_\kappa(t)=A_\mu(t) & = & \frac{q^3tw}{(1-q)^2(1-qt)^2} \cdot C_1 - \frac{q^5t^3w}{(1-q)^2(1-qt)^3} \cdot B_1 \\
& & \mbox{} + \frac{q^5t^3w}{(1-q)^2(1-qt)^4} \cdot A_1 - \frac{q^4t^2w}{(1-q)^2(1-qt)^4} \cdot A(qt) \ .
\end{eqnarray*}
\end{theo}

\noindent \textbf{Proof of Theorem 1.} We shall only prove the formulas for $A_\eta(t)$, $A_\theta(t)$ and $A_\kappa(t)$. The other formulas are proved in a similar way.

Consider this procedure: 

\begin{proc} 
Given a polyomino $P \in S$; nonempty columns $c_1$, $c_2$, $c_3$, $c_5$ and $c_6$, none of which has a hole; a column $c_4$, which may be empty, but may not have a hole. First place the column $c_1$ so that $c_1$ lies two units to the right of the last column of $P$ and so that the top of $c_1$ is aligned with the top of the last column of $P$. Then place the columns $c_2, \ldots, \ c_6$ as shown in Figure 2.
\end{proc}

\begin{figure}
\begin{center}
\includegraphics[width=69mm]{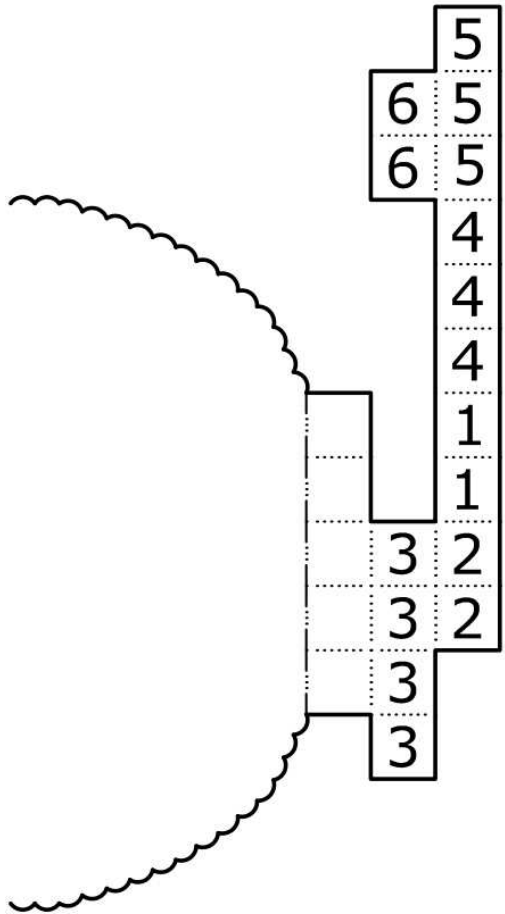}
\caption{An object produced by Procedure 1. For $i=1$ to $6$, the cells of the column $c_i$ are labelled $i$.}
\end{center}
\end{figure}

Procedure 1 does produce all elements of the set $S_\eta$, but also produces some objects which do not lie in $S_\eta$. Indeed, if the column $c_1$ is equally high or higher than the last column of $P$, we get an object in which the second-last and third-last columns do not have any edges in common. That object is obviously not a polyomino. To get rid of those undesired objects, we define a second procedure.

\begin{proc}
Given a polyomino $P \in S$; nonempty columns $c_2$, $c_3$, $c_5$ and $c_6$, none of which has a hole; columns $c_1$ and $c_4$, which may be empty, but may not have a hole. First make a copy of the last column of $P$ and place that copy so that the copy lies two units to the right of the last column of $P$ and so that the top of the copy is aligned with the top of the last column of $P$. Then place the columns $c_1, \ldots, \ c_6$ as shown in Figure 3.
\end{proc}

\begin{figure}
\begin{center}
\includegraphics[width=69mm]{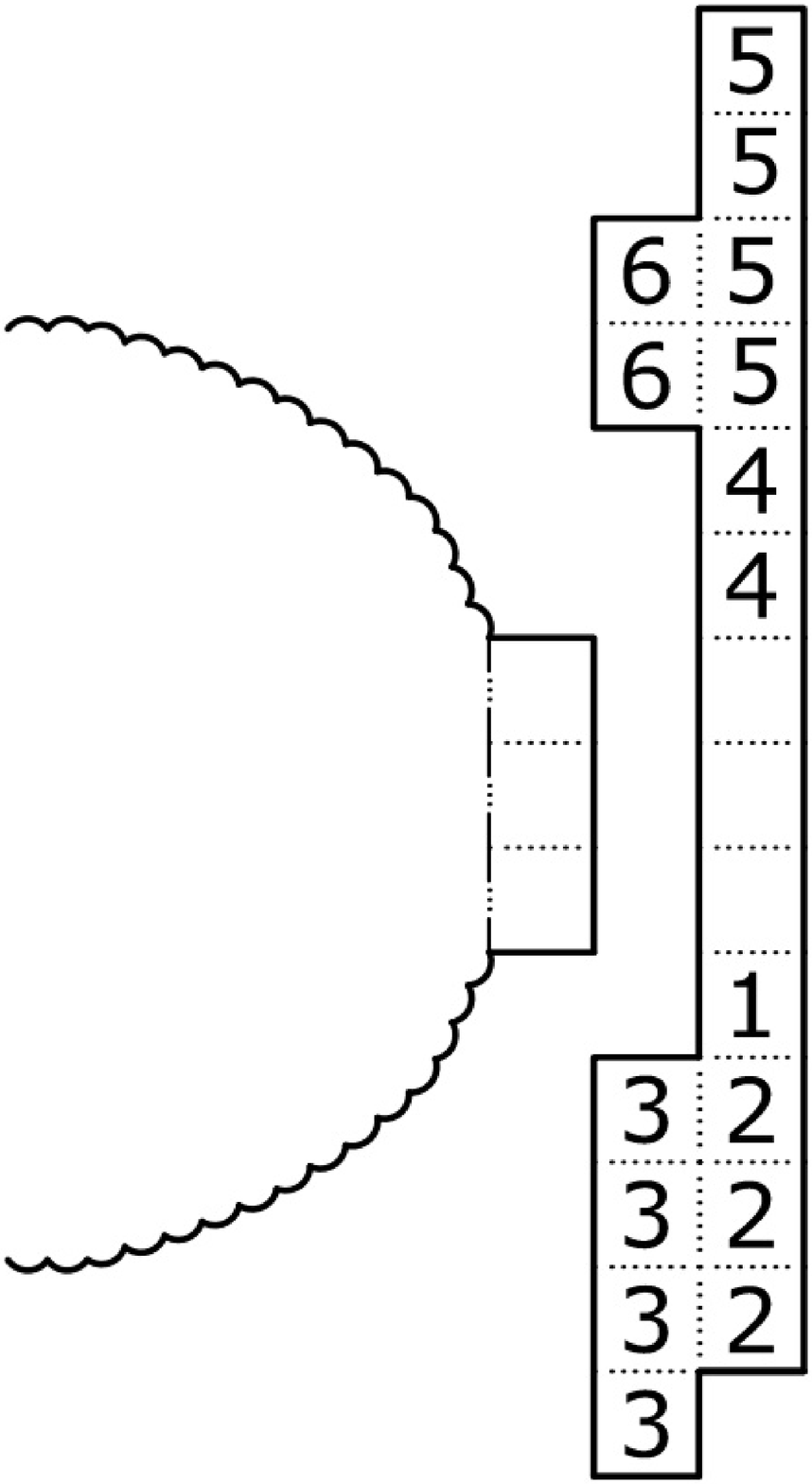}
\caption{An object produced by Procedure 2. For $i=1$ to $6$, the cells of the column $c_i$ are labelled $i$.}
\end{center}
\end{figure}

What Procedure 2 produces is precisely the undesired part of the output of Procedure 1. Now, the generating function for the output of Procedure 1 is 

\begin{equation}
\label{eta1}
A_1 \cdot \left( \frac{q}{1-q} \right)^2 \cdot w \cdot \left( \frac{qt}{1-qt} \right)^3 \cdot \frac{1}{1-qt} 
= \frac{q^5t^3w}{(1-q)^2(1-qt)^4} \cdot A_1.
\end{equation}

In (\ref{eta1}), the contribution of the polyomino $P$ is $A_1$, the contributions of $c_3$ and $c_6$ are $\frac{q}{1-q}$ each, the contributions of $c_1$, $c_2$ and $c_5$ are is $\frac{qt}{1-qt}$ each, and the contribution of $c_4$ is $\frac{1}{1-qt}$. The contributions are multiplied together because Procedure 1 specifies uniquely where each of the pieces is to be placed.

The generating function for the output of Procedure 2 is

\begin{equation}
\label{eta2}
A(qt) \cdot \left( \frac{q}{1-q} \right)^2 \cdot w \cdot \left( \frac{1}{1-qt} \right)^2 \cdot \left( \frac{qt}{1-qt} \right)^2 
= \frac{q^4t^2w}{(1-q)^2(1-qt)^4} \cdot A(qt).
\end{equation}

In (\ref{eta2}), the joint contribution of $P$ and the copy of $P$'s last column is $A(qt)$, the contributions of $c_3$ and $c_6$ are $\frac{q}{1-q}$ each, the contributions of $c_1$ and $c_4$ are $\frac{1}{1-qt}$ each, and the contributions of $c_2$ and $c_5$ are is $\frac{qt}{1-qt}$ each. Again, the contributions are multiplied together because Procedure 2 specifies uniquely where each of the pieces is to be placed.

Subtracting (\ref{eta2}) from (\ref{eta1}) gives

\begin{displaymath}
A_\eta(t) = \frac{q^5t^3w}{(1-q)^2(1-qt)^4} \cdot A_1 - \frac{q^4t^2w}{(1-q)^2(1-qt)^4} \cdot A(qt),
\end{displaymath}

\noindent which was to be demonstrated.

Next, consider this procedure: 

\begin{proc} 
Given a polyomino $P \in S$ and five nonempty columns $c_1, \ldots, \ c_5$, none of which has a hole. First choose a cell (call it $z$) which lies in the last column of $P$ and is not the top cell of that column. Then place the column $c_1$ so that the top cell of $c_1$ becomes the right neighbour of $z$. Then place the columns $c_2, \ldots, \ c_5$ as shown in Figure 4.
\end{proc}

\begin{figure}
\begin{center}
\includegraphics[width=69mm]{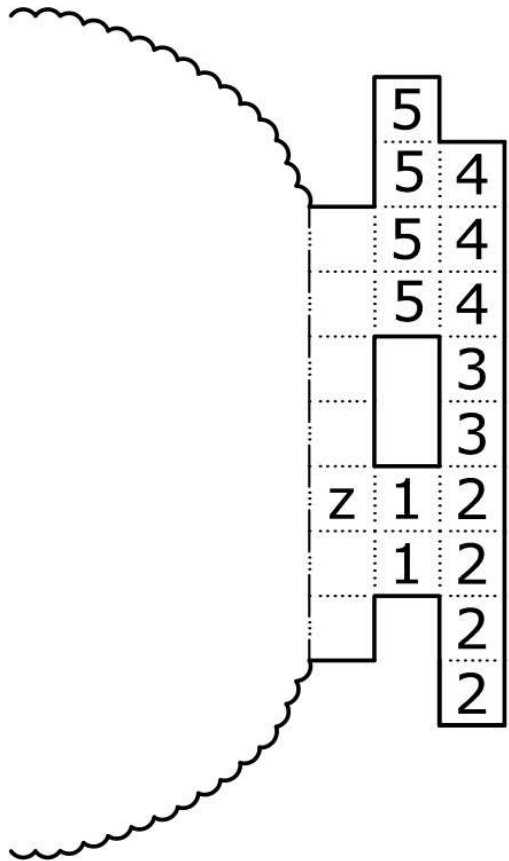}
\caption{An object produced by Procedure 3. For $i=1$ to $5$, the cells of the column $c_i$ are labelled $i$.}
\end{center}
\end{figure}

Procedure 3 produces all elements of the set $S_\theta$, but also produces some objects which do not lie in $S_\theta$. Fortunately, those undesired objects are easy to identify: they are precisely the elements of $S_\eta$, the set whose generating function has just been found. Thus, we shall now compute the generating function for all objects produced by Procedure 3. Then we shall subtract $A_\eta(t)$ and the result will be $A_\theta(t)$.

In Procedure 3, if the last column of $P$ has $n$ cells, the cell $z$ can be chosen in $n-1$ ways. Once $z$ is chosen, the places where $c_1, \ldots, \ c_5$ will go are determined uniquely. The result of Procedure 3 is an object whose second-last column is duplex (the lower and the upper components being $c_1$ and $c_5$, respectively), and the last column is simplex and made up of three parts: $c_2$, $c_3$ and $c_4$. Thus, the objects produced by Procedure 3 have the following generating function:

\begin{eqnarray*}
& \sum_{n=1}^{\infty} (n-1) \cdot \left[ \langle t^n \rangle A(t) \right] \cdot w \cdot \left( \frac{q}{1-q} \right)^2 \cdot \left( \frac{qt}{1-qt} \right)^3 \\
& = \frac{q^5t^3w}{(1-q)^2(1-qt)^3} \cdot \sum_{n=1}^{\infty} (n-1) \cdot \langle t^n \rangle A(t) \\
& = \frac{q^5t^3w}{(1-q)^2(1-qt)^3} \cdot B_1.
\end{eqnarray*}

Consequently,

\begin{eqnarray*}
A_\theta(t) & = & \frac{q^5t^3w}{(1-q)^2(1-qt)^3} \cdot B_1 - A_\eta(t) \\
& = & \frac{q^5t^3w}{(1-q)^2(1-qt)^3} \cdot B_1 - \frac{q^5t^3w}{(1-q)^2(1-qt)^4} \cdot A_1 + \frac{q^4t^2w}{(1-q)^2(1-qt)^4} \cdot A(qt),
\end{eqnarray*}

\noindent as claimed.

We have also promised to prove the formula for $A_\kappa(t)$. So, consider this procedure:
 
\begin{proc} 
Given a polyomino $P \in S$; nonempty columns $c_1$, $c_2$ and $c_3$, none of which has a hole; a column $c_4$, which may be empty, but may not have a hole. First choose two nonadjacent cells which lie in the last column of $P$. Denote the lower chosen cell by $z_1$ and the upper chosen cell by $z_2$. Then place the column $c_1$ so that the top cell of $c_1$ becomes the right neighbour of $z_1$, and place the column $c_2$ so that the bottom cell of $c_2$ becomes the right neighbour of $z_2$. Then place the columns $c_3$ and $c_4$ as shown in Figure 5.
\end{proc}

\begin{figure}
\begin{center}
\includegraphics[width=69mm]{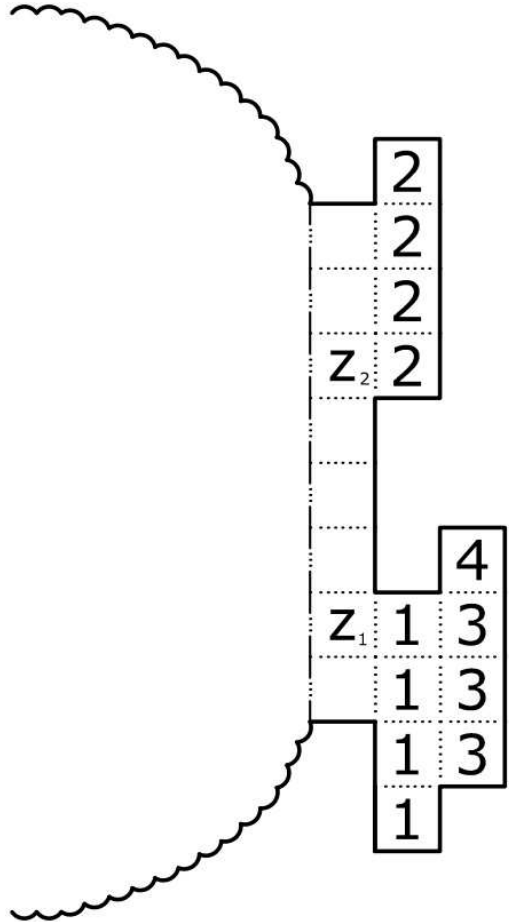}
\caption{An object produced by Procedure 4. For $i=1$ to $4$, the cells of the column $c_i$ are labelled $i$.}
\end{center}
\end{figure}

Procedure 4 produces all elements of the set $S_\kappa$, but also produces some objects which do not lie in $S_\kappa$. Those undesired objects are precisely the elements of $S_\theta$, the set whose generating function has just been found. Thus, we shall now compute the generating function for all objects produced by Procedure 4. Then we shall subtract $A_\theta(t)$ and the result will be $A_\kappa(t)$.

In Procedure 4, if the last column of $P$ has $n$ cells, the cells $z_1$ and $z_2$ can be chosen in ${n-1 \choose 2}=\frac{(n-1)(n-2)}{2}$ ways. Once $z_1$ and $z_2$ are chosen, the places where $c_1, \ldots, \ c_4$ will go are determined uniquely. The result of Procedure 4 is an object whose second-last column is duplex (the lower and the upper components being $c_1$ and $c_2$, respectively), and the last column is simplex and made up of two parts: $c_3$ and $c_4$. Thus, the objects produced by Procedure 4 have the following generating function:

\begin{eqnarray*}
& \sum_{n=1}^{\infty} \frac{(n-1)(n-2)}{2} \cdot \left[ \langle t^n \rangle A(t) \right] \cdot w \cdot \left( \frac{q}{1-q} \right)^2 \cdot \frac{qt}{1-qt}            \cdot \frac{1}{1-qt}\\
& = \frac{q^3tw}{(1-q)^2(1-qt)^2} \cdot \sum_{n=1}^{\infty} \frac{(n-1)(n-2)}{2} \cdot \langle t^n \rangle A(t) \\
& = \frac{q^3tw}{(1-q)^2(1-qt)^2} \cdot C_1.
\end{eqnarray*}

Consequently,

\begin{eqnarray*}
A_\kappa(t) & = & \frac{q^3tw}{(1-q)^2(1-qt)^2} \cdot C_1 - A_\theta(t) \\
& = & \frac{q^3tw}{(1-q)^2(1-qt)^2} \cdot C_1 - \frac{q^5t^3w}{(1-q)^2(1-qt)^3} \cdot B_1 \\ 
& & + \frac{q^5t^3w}{(1-q)^2(1-qt)^4} \cdot A_1 - \frac{q^4t^2w}{(1-q)^2(1-qt)^4} \cdot A(qt), 
\end{eqnarray*}

\noindent as claimed.

\begin{flushright}
$\square$
\end{flushright}

We have $A(t)=A_\alpha(t)+A_\beta(t)+\ldots +A_\mu(t)$, and the expressions for $A_\alpha(t), \ A_\beta(t), \ldots, \ A_\mu(t)$ are given in Theorem 1. Putting these things together, we get a functional equation for $A(t)$. It is convenient to write that functional equation as

\begin{eqnarray}
A(t) & = & \frac{qt}{1-qt} \cdot \left[ 1+B_1+\frac{2q^3w}{(1-q)^3} \cdot C_1 \right] 
+\frac{qt}{(1-qt)^2} \cdot \left[A_1+\frac{2q^2w}{(1-q)^2} \cdot C_1 \right]\nonumber \\
& & \mbox{} +\frac{q^5t^3w}{(1-q)^2(1-qt)^3} \cdot \left(1+\frac{2}{1-q} \cdot A_1 -B_1\right) 
+\frac{3q^5t^3w}{(1-q)^2(1-qt)^4} \cdot A_1 \nonumber \\
& & \mbox{} - \frac{3q^4t^2w}{(1-q)^2(1-qt)^4} \cdot A(qt). \label{eat}
\end{eqnarray}

We solved equation (\ref{eat}) by iteration, as is usually done in the upgraded Temperley method. The iteration ended in

\begin{eqnarray}
A(t) & = & \left\{\sum_{i=1}^{\infty} \frac{(-3)^{i-1}q^{i^2+2i-2}t^{2i-1}w^{i-1}}{(1-q)^{2i-2} \cdot \left[\prod_{k=1}^{i-1}(1-q^kt)\right]^4 \cdot (1-q^it)} \right\} \cdot \left[ 1+B_1+\frac{2q^3w}{(1-q)^3} \cdot C_1 \right] \nonumber \\
& & \mbox{} +\left\{\sum_{i=1}^{\infty} \frac{(-3)^{i-1}q^{i^2+2i-2}t^{2i-1}w^{i-1}}{(1-q)^{2i-2} \cdot \left[\prod_{k=1}^{i-1}(1-q^kt)\right]^4 \cdot (1-q^it)^2} \right\} \cdot \left[A_1+\frac{2q^2w}{(1-q)^2} \cdot C_1 \right] \nonumber \\
& & \mbox{} +\left\{\sum_{i=1}^{\infty} \frac{(-3)^{i-1}q^{i^2+4i}t^{2i+1}w^i}{(1-q)^{2i} \cdot \left[\prod_{k=1}^{i-1}(1-q^kt)\right]^4 \cdot (1-q^it)^3} \right\} \cdot \left(1+\frac{2}{1-q} \cdot A_1 -B_1\right) \nonumber \\
& & \mbox{} -\left\{\sum_{i=1}^{\infty} \frac{(-3)^i q^{i^2+4i}t^{2i+1}w^i}{(1-q)^{2i} \cdot \left[\prod_{k=1}^i(1-q^kt)\right]^4} \right\} \cdot A_1. \label{eatit}
\end{eqnarray}

From equation (\ref{eatit}), we got a system of three linear equations in three unknowns: $A_1$, $B_1$ and $C_1$. The first equation is just the case $t=1$ of (\ref{eatit}). The second equation is obtained by dividing (\ref{eatit}) by $t$, differentiating with respect to $t$ and then setting $t=1$. The third equation is obtained by dividing (\ref{eatit}) by $t$, differentiating twice with respect to $t$, and then setting $t=1$.

Once the linear system is solved, relation (\ref{ekg}) tells us how to obtain the sought-for generating function $G$. The resulting formula is stated in the next section.

\section{A formula for the sought-for generating function} 

\begin{theo} The area and duplex columns generating function for simplex-duplex polyominoes is given by

\begin{equation}
\label{flag}
G(q,w)=\frac{NUM}{DEN},
\end{equation}

\noindent where

\begin{eqnarray*}
NUM & = & (1-q)^4(\tilde{\alpha}+\tilde{\gamma}+2\tilde{\alpha}\tilde{\eta}-2\tilde{\gamma}\tilde{\epsilon})+q^2w(1-q)^2(\tilde{\iota}+\tilde{\lambda}-\tilde{\alpha}\tilde{\kappa}-\tilde{\alpha}\tilde{\mu} \\
& & \mbox{}+\tilde{\beta}\tilde{\iota}+\tilde{\beta}\tilde{\lambda}-\tilde{\gamma}\tilde{\kappa}-\tilde{\gamma}\tilde{\mu}
+\tilde{\delta}\tilde{\iota}+\tilde{\delta}\tilde{\lambda}-2\tilde{\epsilon}\tilde{\lambda}+2\tilde{\eta}\tilde{\iota} +2\tilde{\alpha}\tilde{\zeta}\tilde{\lambda}-2\tilde{\alpha}\tilde{\eta}\tilde{\kappa} \\
& & \mbox{}-2\tilde{\alpha}\tilde{\eta}\tilde{\mu}+2\tilde{\alpha}\tilde{\theta}\tilde{\lambda}-2\tilde{\beta}\tilde{\epsilon}\tilde{\lambda}
+2\tilde{\beta}\tilde{\eta}\tilde{\iota}+2\tilde{\gamma}\tilde{\epsilon}\tilde{\kappa}+2\tilde{\gamma}\tilde{\epsilon}\tilde{\mu}
-2\tilde{\gamma}\tilde{\zeta}\tilde{\iota}-2\tilde{\gamma}\tilde{\theta}\tilde{\iota} \\
& & \mbox{}-2\tilde{\delta}\tilde{\epsilon}\tilde{\lambda}+2\tilde{\delta}\tilde{\eta}\tilde{\iota}) +2q^2w(1-q^2)(\tilde{\alpha}\tilde{\lambda} -\tilde{\gamma}\tilde{\iota}), \\
DEN & = & (1-q)^4(1-\tilde{\beta}+\tilde{\delta}-\tilde{\epsilon}+\tilde{\eta}-\tilde{\alpha}\tilde{\zeta}+\tilde{\alpha}\tilde{\theta}
+\tilde{\beta}\tilde{\epsilon}-\tilde{\beta}\tilde{\eta}+\tilde{\gamma}\tilde{\zeta}
-\tilde{\gamma}\tilde{\theta} \\
& & \mbox{}-\tilde{\delta}\tilde{\epsilon}+\tilde{\delta}\tilde{\eta})-2(1-q)^3(\tilde{\gamma}+\tilde{\alpha}\tilde{\eta}-\tilde{\gamma}\tilde{\epsilon}) \\
& & \mbox{}-2q^2w(1-q)^2(\tilde{\kappa}-\tilde{\beta}\tilde{\mu}+\tilde{\delta}\tilde{\kappa}-\tilde{\epsilon}\tilde{\kappa}+\tilde{\zeta}\tilde{\iota} 
-\tilde{\zeta}\tilde{\lambda}+\tilde{\eta}\tilde{\kappa}-\tilde{\alpha}\tilde{\zeta}\tilde{\mu}+\tilde{\alpha}\tilde{\theta}\tilde{\kappa} \\
& & \mbox{}+\tilde{\beta}\tilde{\epsilon}\tilde{\mu}-\tilde{\beta}\tilde{\eta}\tilde{\mu}-\tilde{\beta}\tilde{\theta}\tilde{\iota}
+\tilde{\beta}\tilde{\theta}\tilde{\lambda}+\tilde{\gamma}\tilde{\zeta}\tilde{\mu}-\tilde{\gamma}\tilde{\theta}\tilde{\kappa}
-\tilde{\delta}\tilde{\epsilon}\tilde{\kappa}+\tilde{\delta}\tilde{\zeta}\tilde{\iota}-\tilde{\delta}\tilde{\zeta}\tilde{\lambda} +\tilde{\delta}\tilde{\eta}\tilde{\kappa}) \\
& & \mbox{}-4q^2w(1-q)(\tilde{\beta}\tilde{\lambda}-\tilde{\gamma}\tilde{\kappa}+\tilde{\alpha}\tilde{\zeta}\tilde{\lambda}   -\tilde{\alpha}\tilde{\eta}\tilde{\kappa}-\tilde{\beta}\tilde{\epsilon}\tilde{\lambda}+\tilde{\beta}\tilde{\eta}\tilde{\iota} +\tilde{\gamma}\tilde{\epsilon}\tilde{\kappa}-\tilde{\gamma}\tilde{\zeta}\tilde{\iota}) \\
& & \mbox{}-2q^3w(1-q)(\tilde{\iota}+\tilde{\alpha}\tilde{\kappa}-\tilde{\alpha}\tilde{\mu}-\tilde{\beta}\tilde{\iota}+\tilde{\delta}\tilde{\iota}
-\tilde{\epsilon}\tilde{\lambda}+\tilde{\eta}\tilde{\iota}-\tilde{\alpha}\tilde{\zeta}\tilde{\lambda}+\tilde{\alpha}\tilde{\eta}\tilde{\kappa} \\
& & \mbox{}-\tilde{\alpha}\tilde{\eta}\tilde{\mu}+\tilde{\alpha}\tilde{\theta}\tilde{\lambda}+\tilde{\beta}\tilde{\epsilon}\tilde{\lambda}
-\tilde{\beta}\tilde{\eta}\tilde{\iota}-\tilde{\gamma}\tilde{\epsilon}\tilde{\kappa}+\tilde{\gamma}\tilde{\epsilon}\tilde{\mu}
+\tilde{\gamma}\tilde{\zeta}\tilde{\iota}-\tilde{\gamma}\tilde{\theta}\tilde{\iota}-\tilde{\delta}\tilde{\epsilon}\tilde{\lambda}
+\tilde{\delta}\tilde{\eta}\tilde{\iota}) \\
& & -4q^3w(\tilde{\alpha}\tilde{\lambda}-\tilde{\gamma}\tilde{\iota}),
\end{eqnarray*}

\begin{eqnarray*}
\tilde{\beta} & = & \sum_{i=1}^{\infty} \frac{(-3)^{i-1}q^{i^2+2i-2}w^{i-1}}{(1-q)^{2i-2} \cdot \left[\prod_{k=1}^{i-1}(1-q^k)\right]^4 \cdot 
(1-q^i)^{{\scriptstyle{\overline{2}}}\atop}} \ , \\
\tilde{\gamma} & = & \sum_{i=1}^{\infty} \frac{(-3)^{i-1}q^{i^2+4i}w^i}{(1-q)^{2i} \cdot \left[\prod_{k=1}^{i-1}(1-q^k)\right]^4 \cdot 
(1-q^i)^{{\scriptstyle{\overline{3}}}\atop}} \ , \\
\tilde{\zeta} & = & \sum_{i=1}^{\infty} \frac{(-3)^{i-1}q^{i^2+2i-2}w^{i-1}}{(1-q)^{2i-2} \cdot \left[\prod_{k=1}^{i-1}(1-q^k)\right]^4 \cdot 
(1-q^i)^{{\scriptstyle{\overline{2}}}\atop}} \\ 
& & \cdot \left(2i-2+4 \cdot \sum_{k=1}^{i-1}\frac{q^k}{1-q^k} + \frac{\overline{2}q^i}{1-q^i} \right) , \\
\tilde{\eta} & = & \sum_{i=1}^{\infty} \frac{(-3)^{i-1}q^{i^2+4i}w^i}{(1-q)^{2i} \cdot \left[\prod_{k=1}^{i-1}(1-q^k)\right]^4 \cdot 
(1-q^i)^{{\scriptstyle{\overline{3}}}\atop}} \\ 
& & \cdot \left(2i+4 \cdot \sum_{k=1}^{i-1}\frac{q^k}{1-q^k} + \frac{\overline{3}q^i}{1-q^i} \right) , \\
\tilde{\kappa} & = & \frac{1}{2} \cdot \sum_{i=1}^{\infty} \frac{(-3)^{i-1}q^{i^2+2i-2}w^{i-1}}{(1-q)^{2i-2} \cdot \left[\prod_{k=1}^{i-1}(1-q^k)\right]^4 \cdot 
(1-q^i)^{{\scriptstyle{\overline{2}}}\atop}} \\ 
& & \cdot \left[\left(2i-2+4 \cdot \sum_{k=1}^{i-1}\frac{q^k}{1-q^k} + \frac{\overline{2}q^i}{1-q^i} \right)^2 \right. \\
& & \left. -2i+2+4 \cdot \sum_{k=1}^{i-1}\frac{q^{2k}}{(1-q^k)^2} + \frac{\overline{2}q^{2i}}{(1-q^i)^2}\right] , \\
\tilde{\lambda} & = & \frac{1}{2} \cdot \sum_{i=1}^{\infty} \frac{(-3)^{i-1}q^{i^2+4i}w^i}{(1-q)^{2i} \cdot \left[\prod_{k=1}^{i-1}(1-q^k)\right]^4 \cdot 
(1-q^i)^{{\scriptstyle{\overline{3}}}\atop}} \\ 
& & \cdot \left[\left(2i+4 \cdot \sum_{k=1}^{i-1}\frac{q^k}{1-q^k} + \frac{\overline{3}q^i}{1-q^i} \right)^2 \right. \\
& & \left. -2i+4 \cdot \sum_{k=1}^{i-1}\frac{q^{2k}}{(1-q^k)^2} + \frac{\overline{3}q^{2i}}{(1-q^i)^2} \right] .
\end{eqnarray*}

In the formula for $\tilde{\beta}$, one of the numbers has an overline. (Having an overline does not affect the value of a number.) To obtain the formula for $\tilde{\alpha}$ from the formula for $\tilde{\beta}$, it is enough to replace the said $\overline{2}$ by $1$. 

In the formula for $\tilde{\gamma}$, there is one $\overline{3}$. To obtain the formula for $\tilde{\delta}$ from the formula for $\tilde{\gamma}$, it is enough to replace the $\overline{3}$ by $4$ and change $(-3)^{i-1}$ into $(-3)^i$.

In the formula for $\tilde{\zeta}$ (resp. $\tilde{\kappa}$), there are two (resp. three) $\overline{2}$'s. To obtain the formula for $\tilde{\epsilon}$ from the formula for $\tilde{\zeta}$, and also to obtain the formula for $\tilde{\iota}$ from the formula for $\tilde{\kappa}$, it is enough to replace each of the $\overline{2}$'s by $1$.

In the formula for $\tilde{\eta}$ (resp. $\tilde{\lambda}$), there are two (resp. three) $\overline{3}$'s. To obtain the formula for $\tilde{\theta}$ from the formula for $\tilde{\eta}$, and also to obtain the formula for $\tilde{\mu}$ from the formula for $\tilde{\lambda}$, it is enough to replace each of the $\overline{3}$'s by $4$ and change $(-3)^{i-1}$ into $(-3)^i$.

\end{theo}

By setting $w=0$, from Theorem 2 we obtain a widely known result, discovered independently by Temperley \cite{Temperley} and P\'{o}lya \cite{Polya} in the early days of polyomino enumeration:

\begin{coro}
The area generating function for column-convex polyominoes is given by
\begin{displaymath}
G(q,0)=\frac{q(1-q)^3}{1-5q+7q^2-4q^3} \ .
\end{displaymath}
\end{coro}

\noindent \textbf{Proof.} A short calculation.
\begin{flushright}
$\square$
\end{flushright}

\section{A bit of asymptotic analysis}

Formula (\ref{flag}) is very complicated, but still allows us to quickly compute Taylor polynomials of any reasonable degree. We actually computed the Taylor polynomial of degree $320$ for the function $G(q,1)$. It turned out that 

\begin{eqnarray*}
G(q,1) & = & q+2q^2+6q^3+19q^4+63q^5+216q^6+758q^7+2693q^8 \\
& & \mbox{}+9608q^9+34269q^{10}+121946q^{11}+432701q^{12}+\ldots \ .
\end{eqnarray*}

Then, for $n=2,\ 3,\ldots,\ 320$, we divided the coefficient $\langle q^n \rangle G(q,1)$ by the coefficient $\langle q^{n-1} \rangle G(q,1)$. This quotient gradually stabilizes, so that for $n=91,\ 92,\ldots,\ 320$ we have 

\begin{displaymath}
\left[\frac{\langle q^n \rangle G(q,1)}{\langle q^{n-1} \rangle G(q,1)}\right]_{rounded \ to \ 12 \ decimal \ places}=3.522019812882 \ .
\end{displaymath}

Next, for $n=1,\ 2,\ldots,\ 320$, we divided $\langle q^n \rangle G(q,1)$ by $3.5220198128815848301^n$. Once again, the quotient gradually
stabilizes: for $n=78,\ 79,\ldots,\ 320$ we have $\frac{\langle q^n \rangle G(q,1)}{3.5220198128815848301^n}=0.119442870405$. So, there is plenty of evidence 
that the coefficient $\langle q^n \rangle G(q,1)$ (\textit{i.e.}, the number of $n$-celled simplex-duplex polyominoes) has the asymptotic behaviour 
$\langle q^n \rangle G(q,1) \sim 0.119442870405 \cdot 3.522019812882^n$. The dominant singularity of $G(q,1)$ is a simple pole, located at $0.283927988236$.

For comparison, $\langle q^n \rangle G(q,0)$ (\textit{i.e.}, the number of $n$-celled column-convex polyominoes) has the asymptotic behaviour
$\langle q^n \rangle G(q,1) \sim 0.180915501882 \cdot 3.205569430401^n$. The dominant singularity of $G(q,0)$ is a simple pole, located at $0.311957055279$.

\section{Further work}

Our next goal is to find the area generating function for simplex-duplex polyominoes with hexagonal cells. That should not be difficult because we now have a method which works with square cells. \textit{Mutatis mutandis}, such methods usually also work when cells are hexagons.

We are sceptical about the feasibility of an exact \textit{perimeter} enumeration of simplex-duplex polyominoes. Anyway, it would be worthwhile to effectively do the computations and thus advance from scepticism to certainty. Further, it might be possible (although we are sceptical again) to compute the area generating function of \textit{simplex-multiplex polyominoes}. In a simplex-multiplex polyomino, the columns may have any number of connected components, but it is not allowed that two adjacent columns have two or more connected components each.

On the other hand, it is probably possible to compute the area generating function of \textit{simplex-duplex}$_2$ polyominoes. Here, by a 
simplex-duplex$_2$ polyomino we mean a column-duplex polyomino in which runs of two consecutive two-component columns are allowed, but it is not allowed that three consecutive columns have two connected components each. We feel, however, that the game is not worth the candle. Namely, simplex-duplex polyominoes already have a very complicated area generating function, and simplex-duplex$_2$ polyominoes will definitely have a still more complicated area generating function.

\end{document}